\documentclass[11pt]{article}
\usepackage{a4}
\usepackage{latexsym}
\usepackage{amsfonts}
\usepackage{amssymb}
\usepackage{epic}
\usepackage{xypic}
\usepackage[dvips]{graphicx}
\usepackage{amscd}
\usepackage{amsmath}
\usepackage[noxcolor]{pstricks}
\usepackage[normalem]{ulem}
\newtheorem{theorem}{Theorem}
\newtheorem{lemma}[theorem]{Lemma}
\newtheorem{proposition}[theorem]{Proposition}
\newtheorem{corollary}[theorem]{Corollary}
\newtheorem{definition}
{Definition}
\newtheorem{remark}{Remark\/}
\newtheorem{example}{Example\/}
\newtheorem{conjecture}
{Conjecture\/}

\newcommand{\C}{{\mathbb C}}

\begin{document}
\title{\textbf{The monodromy conjecture for plane meromorphic germs}}
\author{Manuel Gonz\'alez Villa and Ann Lemahieu}  \date{}
\maketitle {\footnotesize \emph{\textbf{Abstract.---}
A notion of Milnor fibration  for meromorphic functions and the corresponding concepts of  monodromy and monodromy zeta function have been introduced in {\rm \cite{meromorphic}}. In this article we define the topological zeta function for meromorphic germs and we study its poles in the plane case. We show that the poles do not behave as in the holomorphic case but still do satisfy a generalization of the monodromy conjecture.}}

${}$
\begin{center}
\textsc{0. Introduction}
\end{center}
${}$ \\
In the last decades there has been an increasing interest in
zeta functions associated to a germ of holomorphic function, such as the Igusa zeta function, the topological zeta function and the motivic zeta function. These zeta functions are rational and can be described in terms of an embedded resolution of the germ. Each irreducible component of the total transform of an embedded resolution of the germ yields a `candidate pole' of these zeta functions. However, many of these candidate poles are canceled.  In general no geometric criterion is known to sort out the poles from a set of candidate poles. The poles are completely characterised in dimension two. In \cite[Theorem 4.3]{Veys} Veys gives a criterion to deduce easily the poles from the embedded resolution graph.

The poles of these zeta functions also show up in the monodromy conjecture. The monodromy conjecture predicts that if $s_0$ is a pole of the topological zeta function (or the Igusa or the motivic zeta function) associated to a holomorphic germ $f: (\C^{n+1},0) \rightarrow (\C,0)$, then $e^{2\pi i s_0}$ is an eigenvalue of monodromy at some point of $f^{-1}\{0\}$ in a neighbourhood of the origin.
This conjecture has been proven for plane curves by Loeser in \cite{loeser1}. An alternative proof is given by Rodrigues in \cite{bart}.
The conjecture is still open in arbitrary dimension. Particular cases are proven in \cite{loeser2}, \cite{luengo}, \cite{luengo2}, \cite{monodconjtoric} and \cite{liseann}.

These poles are rational numbers of the form $-\nu_i/N_i, i \in S,$ where the numbers $\nu_i-1$ and $N_i$ are given by the multiplicities of the divisors $\pi^\ast dx_1 \wedge \cdots \wedge dx_{n+1}$ and  $\pi^\ast f$ along $E_i$ for an embedded resolution $\pi: X \rightarrow \mathbb{C}^{n+1}$ of $f$
(see Section 1).  N\'emethi and Veys propose an original approach to the conjecture (at least for dimension 2)  in \cite{Veys4}, \cite{wimnemethi1} and \cite{wimnemethi2}. By considering another differential form than the standard form a different set of numbers $(\nu_i)_{i \in S}$ is generated. They investigate under which conditions the corresponding new set of poles still gives rise to eigenvalues of monodromy and which eigenvalues can be found in this way.

In this article we take a complementary approach. In our case the values  $(\nu_i)_{i \in S}$ are the standard ones but we vary the set $(N_i)_{i \in S}$. In concrete, we consider the case of plane meromorphic germs $f=\frac{P}{Q}$. The series of numbers $(N_i =N^P_i - N^Q_i)_{i \in S}$ are now more general than in the case of plane curves.
The notions of Milnor fibration and monodromy for germs of meromorphic functions have been defined in \cite{meromorphic}. These definitions were motivated by the classification of germs of fractions due to Arnold in \cite{arnold}.
We define a topological zeta function for meromorphic germs and then, the context for the monodromy conjecture being prepared, we show that the monodromy conjecture also holds for plane meromorphic germs.
Very recently, Raibaut defined a motivic Milnor fibre and motivic zeta function for rational functions (see \cite{raibaut, raibautb}).

The outline of this paper is as follows. In the first section we recall the notions of Milnor fibre for meromorphic germs and zeta function of monodromy, as they were introduced by Gusein-Zade, Luengo and Melle-Hern\'andez in \cite{meromorphic}. In the second section we define the topological zeta function in the meromorphic
context and we study its poles in the plane case. We show that now just only one implication in the criterion of Veys holds.
In Section 3 we propose the monodromy conjecture for meromorphic germs and we prove the conjecture in the plane case.

 \begin{center}
\textsc{1. Monodromy for meromorphic germs}
\end{center}

\begin{definition}
A \emph{germ of a meromorphic function} on $(\mathbb{C}^{n+1},0)$ is a fraction $f=\frac{P}{Q}$, where $P$ and $Q$ are germs of holomorphic functions $P,Q: (\mathbb{C}^{n+1},0) \rightarrow (\mathbb{C},0)$.
Two meromorphic functions $f=P/Q$ and $f'=P'/Q'$ are said to be equal if there exists a germ of a holomorphic function $u: (\mathbb{C}^{n+1},0) \rightarrow (\mathbb{C},0)$ with $u(0) \neq 0$ such that $P'=Pu$ and $Q'=Qu$.
\end{definition}

\begin{definition}
Let $X$ be an $(n+1)$-dimensional smooth analytic manifold, $\mathcal{U}$ a neighbourhood of $0 \in \mathbb{C}^{n+1}$ and $\pi: X \rightarrow \mathcal{U}$ a proper analytic map which is an isomorphism outside a proper analytic subspace in $\mathcal{U}$. We say that $\pi$ is an \emph{embedded resolution of the germ of the meromorphic function} $f=\frac{P}{Q}$ if
\begin{enumerate}
\item the total transform $\pi^{-1}(H)$ of the hypersurface $H=\{P_{|\mathcal{U}}=0\} \cup \{Q_{|\mathcal{U}}=0\}$ is a normal crossing divisor at each point of $X$;
\item the lifting $\tilde{f}=f \circ \pi=\frac{P \circ \pi}{Q \circ \pi}$ defines a holomorphic map $\tilde{f}: X \rightarrow \mathbb{P}^1$.
\end{enumerate}
\end{definition}

\begin{remark}
\emph{One can obtain an embedded resolution $\pi$ for $\frac{P}{Q}$ from an embedded resolution $\pi'$ for $PQ$ by blowing up along the intersections of irreducible components of $\pi^{-1}\{PQ=0\}$ until the irreducible components of the strict transform in $\pi^{-1}\{P=0\}$ and $\pi^{-1}\{Q=0\}$ are separated by a \emph{dicritical component}, i.e.\ an exceptional divisor $E$ of $\pi$ for which $\tilde{f}_{|E}: E \rightarrow \mathbb{P}^1$ is a surjective map.}
\end{remark}

\begin{example}
\emph{We compute an embedded resolution $\pi$ for the germ $f=P/Q=((y^2-x^3)^2-xy^5)/(x-y)$ and represent it by its \emph{dual resolution graph}, i.e.\ to each irreducible component of $\pi^{-1}(H)$ we associate a vertex and we join two vertices if the corresponding components have a non-empty intersection. The vertices corresponding to the components of the strict transform of $H$ are represented by arrowheads. An arrowhead  $\rightarrow$ (resp. $\twoheadrightarrow$) corresponds to an irreducible component of the strict transform of $P$ (resp. $Q$). The following figure shows the dual resolution graph of $f$.
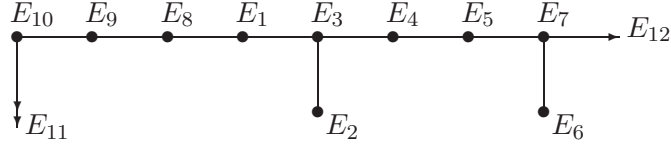
\begin{center}
\begin{picture}(90,30)(0,0)
\put(9,19){$\bullet$} \put(19,19){$\bullet$} \put(29,19){$\bullet$} \put(39,19){$\bullet$} \put(49,19){$\bullet$} \put(59,19){$\bullet$} \put(69,19){$\bullet$} \put(79,19){$\bullet$}
\put(10,20){\line(1,0){70}}
\put(10,20){\vector(0,-1){10}}\put(10,20){\vector(0,-1){12}}
\put(50,20){\line(0,-1){10}}\put(49,9){$\bullet$}
\put(80,20){\line(0,-1){10}}\put(79,9){$\bullet$} \put(80,20){\vector(1,0){10}}
\put(9,22){$E_{10}$}
\put(19,22){$E_{9}$}\put(29,22){$E_{8}$}\put(39,22){$E_{1}$}\put(49,22){$E_{3}$}\put(51,7){$E_{2}$}\put(59,22){$E_{4}$}
\put(69,22){$E_{5}$}\put(79,22){$E_{7}$}\put(81,7){$E_{6}$}\put(11,7){$E_{11}$}\put(91,20){$E_{12}$}
\end{picture}
\end{center}
\vspace*{-0.7cm}
The exceptional components $E_1, \cdots, E_7$ are created to resolve $\{PQ=0\}$.
We then perform some additional blow ups until we get the dicritical component $E_{10}$ that
separates the components according to the value they take when applying $f \circ \pi$. In this example
$E_{11}$ is send to $\infty \in \mathbb{P}^1$ and the other components different from the dicritical component $E_{10}$ are mapped to $0 \in \mathbb{P}^1$.}\hfill $\square$
\end{example}
In \cite{meromorphic2}, Gusein-Zade, Luengo and Melle-Hern\'andez introduced the notion
of $a$-Milnor fibre for meromorphic germs. Let $f=P/Q$ be a meromorphic germ on $(\C^{n+1},0)$.
Let $B_{\epsilon}$ be the closed ball of radius $\epsilon$ centered at $0 \in \C^{n+1}$ with $\epsilon$ sufficiently small such that $P$ and $Q$ are defined on $B_{\epsilon}$ and the sphere $\partial B_{\epsilon}$ intersects transversally $\{P=0\}$ and $\{Q=0\}$. Then $f : B_{\epsilon} \setminus \{P=Q=0\} \rightarrow \mathbb{P}^1$ is a locally trivial $C^{\infty}$ fibration over a small enough punctured disc with center $a \in \mathbb{P}^1$.
\begin{definition}
The \emph{$a$-Milnor fibre of the meromorphic germ $f$ at $0 \in \C^{n+1}$} is
\[\mathcal{M}^a_{f} = \{z \in B_{\epsilon} \mid f(z)=\frac{P(z)}{Q(z)}=c\},\]
with $c \in \C$ and $||c-a||$ small enough and different from 0.
\end{definition}
They also introduced the corresponding $a$-monodromy zeta function.
Consider a counterclockwise going loop around the value $a$
and let $h^a_{f}$ 
be the induced diffeomorphism of the $a$-Milnor fibre $\mathcal{M}^a_{f}$. The \emph{$a$-monodromy transformation $h^a_{f^*}$}
at the origin is then the induced automorphism of the homology groups of the $a$-Milnor fibre and its zeta function is defined as
\begin{eqnarray*}
\zeta^a_{f}(t) & = & \prod_{q \geq 0} \left(\textnormal{det}({Id - th^a_{f^*}}_{|H_q(\mathcal{M}^{a}_{f},\C)}) \right)^{(-1)^q}. 
\end{eqnarray*}
Notice that $\infty$-Milnor fibre of $f$ is the $0$-Milnor fibre of $1/f$ (and $\zeta^{\infty}_{f}=\zeta^0_{1/f}$) and that the $a$-Milnor fibre of $f$, for $a \in \mathbb{P}^1 \setminus \{0,\infty\}$, is the $0$-Milnor fibre of $f-a$ (and then $\zeta^a_{f}=\zeta^0_{f-a}$). Therefore in the rest of the paper we permit us to think at $a=0$ and we will not explicitly denote the value $a=0$.

In \cite{meromorphic}, the monodromy zeta function is expressed  in terms of an embedded resolution of $f$. We first  introduce some useful notation. Let $\pi: X \rightarrow \mathcal{U}$ be an embedded resolution of the germ of meromorphic function $f=\frac{P}{Q}$ and let $(E_i)_{i \in S}$ be the irreducible components of the total transform $\pi^{-1}(H)$.
We will denote $E_i^{\circ}:= E_i \setminus
(\cup_{j \in S \setminus \{i\}} E_j)$, for $i \in S$.
For $i \in S$, let $N_i^P$ (resp. $N_i^Q$) be the multiplicity of $P \circ \pi$ (resp. $Q \circ \pi$) along $E_i$ and $\nu_i-1$ the multiplicity of $\pi^*(dx_1 \wedge \cdots \wedge dx_{n+1})$ along $E_i$. The data $(N_i^P,N_i^Q,\nu_i)$ are called the \emph{numerical data} corresponding to the component $E_i$. Sometimes we just write $(N_i,\nu_i)$, with $N_i:=N_i^P-N_i^Q$. Finally, denote by  $S_{0}$ the set $\{i \in S \mid N_i^P > N_i^Q\}$.

Remark that the irreducible components of the total transform that are send to $0$ (resp. to $\infty$) under $f \circ \pi$ are then exactly those $E_i$ for which $N_i^P > N_i^Q$ (resp. $N_i^P < N_i^Q$). The dicritical components verify that $N_i^P=N_i^Q$. However $N_i^P=N_i^Q$ does not imply that $E$ is dicritical.

We denote the topological Euler characteristic by $\chi(\cdot)$.
\begin{theorem}{\rm \cite[Theorem 1]{meromorphic}}
Let $\pi: X \rightarrow \mathcal{U}$ be an isomorphism outside the hypersurface $\{P_{|\mathcal{U}}=0\} \cup \{Q_{|\mathcal{U}}=0\}$. Then
\begin{eqnarray*}
\zeta_{f}(t) & = & \prod_{i \in S_0} \left(1-t^{N_i^P - N_i^Q}\right)^{\chi(E_i^{\circ} \cap \pi^{-1}\{0\})}, 
\end{eqnarray*}
\end{theorem}
\textbf{Example 1 Continued}
The numerical data corresponding to the embedded resolution $\pi$ for the germ $f=P/Q=((y^2-x^3)^2-xy^5)/(x-y)$ of Example 1 are written in the following table.
\\ ${}$ \begin{center} \begin{tabular}{|c|c|c|c|c|c|c|c|c|c|c|c|c|}
  \hline
  $i=$ & $1$ & $2$ & $3$ & $4$ & $5$ & $6$ & $7$ & $8$ & $9$ & $10$ & $11$ & $12$ \\
   \hline
  $N_i^P$ & 4 & 6 & 12 & 14 & 16 & 17 & 34 & 4 & 4 & 4 & 0 & 1 \\
  $N_i^Q$ & 1 & 1 & 2 & 2 & 2 & 2 & 4 & 2 & 3 & 4 & 1 & 0\\
  $\nu_i$ & 2 & 3 & 5 & 6 & 7 & 8 & 15 & 3 & 4 & 5 & 1 & 1\\
  \hline
\end{tabular}
\end{center}
We thus find
\[\zeta_{f}(t) = \frac{(1-t^5)(1-t^{15})}{(1-t^{10})(1-t^{30})}.\] \hfill $\square$

\begin{center}
\textsc{2. Topological zeta function for meromorphic functions}
\end{center}
${}$\\
In this section we extend the definition of the local topological zeta function for holomorphic functions to meromorphic functions. 
Let $f=\frac{P}{Q}$ be a meromorphic function on $(\C^{n+1},0)$ and let $\pi: X \rightarrow \mathcal{U}$ be an embedded resolution of $f$ as in Section 1. We denote $E_I
:= \cap_{i \in I} E_i$ and $E_I^{\circ}:= E_I \setminus (\cup_{j
\notin I} E_j)$.
\begin{definition}
The \emph{0-local topological zeta function associated to $f$ at $0 \in \C^{n+1}$} is the rational function
\[Z_{top,f}(s)=\sum_{\substack{I \subset S \\ I \cap S_{0} \neq \emptyset}} \chi(E_I^{\circ}
\cap \pi^{-1}\{ 0 \})\prod_{i \in I} \frac{1}{(N_i^P-N_i^Q)s+ \nu_i}.\]
\end{definition}
As in the holomorphic case, the global version is defined by taking $\chi(E_I^{\circ})$ instead of $\chi(E_I^{\circ} \cap \pi^{-1}\{ 0 \})$. Checking that this definition does not depend on the chosen embedded resolution for $f$ goes analogously as in the holomorphic case. \\
The definition of the topological zeta function yields a set of candidate poles
\[\left\{-\frac{\nu_i}{N_i^P-N_i^Q} \mid i \in S_0\right\}\]
that are negative rational numbers, depending on the chosen embedded resolution. Notice that in particular the dicritical components do not give rise to candidate poles and neither components of the strict transform of $Q$ do.
\\ \\ \noindent \textbf{Example 1 Continued Continued.}
The local topological zeta function of $f=\frac{(y^2-x^3)^2-xy^5}{x-y}$ at $0$ is
\begin{small}\begin{eqnarray*}
 Z_{top,f}(s) & = & -\frac{1}{10s+5}-\frac{1}{30s+15}+\frac{1}{5(s+4)}+\frac{1}{(s+4)(2s+3)}+\frac{1}{(2s+3)(3s+2)} \\
& + & \frac{1}{(10s+5)(3s+2)}+\frac{1}{(10s+5)(5s+3)}+\frac{1}{(10s+5)(12s+6)} \\ & + & \frac{1}{(12s+6)(14s+7)}
+ \frac{1}{(14s+7)(30s+15)}+\frac{1}{(30s+15)(15s+8)} \\ & + & \frac{1}{(30s+15)(s+1)}  +  \frac{1}{15s+8}+\frac{1}{5s+3}
\\ & = & \frac{20s^2+33s+12.}{15(s+1)(2s+1)^2}
\end{eqnarray*}\end{small}
\hfill $\square $
\\ \\ In the above example we see that most candidate poles are canceled. This is a general phenomenon and it appears already in the holomorphic case. For dimension 2 there is a complete geometric criterion to determine the poles of the local topological zeta function (or the other local zeta functions) associated to a holomorphic function. The general case is however poorly understood. We recall this criterion.
\begin{theorem}{\rm \cite{Veys}} \label{thmVeys}
Let $f \in \C[x,y]$ and let $\pi$ be the minimal embedded resolution of $f$.
We have that $s_0$ is a pole of $Z_{top,f}(s)$ if and only if $s_0 = -\nu_i/N_i$ for
some exceptional curve $E_i$ intersecting at least three times other components, or
$s_0 = -1/N_j$ for some irreducible component $E_j$ of the strict transform of $f^{-1}\{0\}$.
\end{theorem}
We will now study the poles in the case of plane meromorphic germs.
The next example shows that Theorem \ref{thmVeys} does not hold anymore for plane meromorphic functions.

\begin{example}
\emph{Let $f=P/Q$ with
\normalsize\[P(x,y)=\prod_{i=1}^3 \prod_{j=1}^5((y+a_ix)-b_jx^2)P'(x,y) \quad \textnormal{and}\]
\scriptsize{\[P'(x,y)=\prod_{k=1}^4\left(\prod_{l=1}^2((y+c_kx)-x^2-c_kd_lx^2-d_lxy) \prod_{m=1}^3((y+c_kx)-2x^2-c_kf_mx^2-f_mxy)\right),\]}
\normalsize\[Q(x,y)=\prod_{n=1}^{35}(y-g_nx),\]
with $a_i,c_k,g_n \in \C$ all different, all the $b_j \in \C$ different, the $d_l \in \C$ different and the $f_m \in \C$ different.
The dual resolution graph is as follows
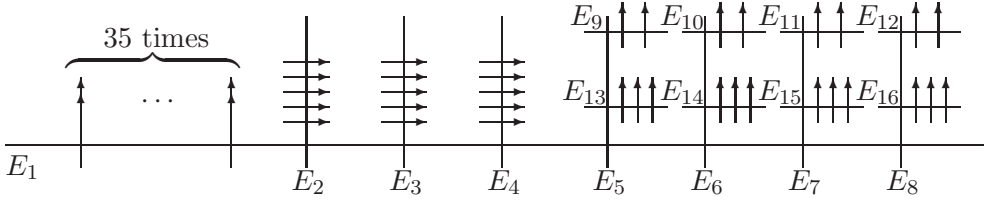
\begin{center}
\begin{picture}(130,30)(0,10)
\put(0,20){\line(1,0){130}} \put(10,17){\vector(0,1){10}}\put(10,17){\vector(0,1){12}}
\put(18,25){$\cdots$}  \put(30,17){\vector(0,1){10}} \put(30,17){\vector(0,1){12}}
\put(8,30){$\overbrace{\qquad \qquad \qquad }$}\put(13,33){35 times}
\put(40,17){\line(0,1){20}} \put(53,17){\line(0,1){20}} \put(66,17){\line(0,1){20}}
\put(37,25){\vector(1,0){6}} \put(37,27){\vector(1,0){6}}\put(37,29){\vector(1,0){6}}
\put(37,31){\vector(1,0){6}} \put(37,23){\vector(1,0){6}}
\put(50,25){\vector(1,0){6}} \put(50,27){\vector(1,0){6}}\put(50,29){\vector(1,0){6}}
\put(50,31){\vector(1,0){6}} \put(50,23){\vector(1,0){6}}
\put(63,25){\vector(1,0){6}} \put(63,27){\vector(1,0){6}}\put(63,29){\vector(1,0){6}}
\put(63,31){\vector(1,0){6}} \put(63,23){\vector(1,0){6}}
\put(80,17){\line(0,1){20}} \put(93,17){\line(0,1){20}} \put(106,17){\line(0,1){20}}
\put(119,17){\line(0,1){20}}
\put(77,25){\line(1,0){11}}\put(77,35){\line(1,0){11}}
\put(90,25){\line(1,0){11}}\put(90,35){\line(1,0){11}}
\put(103,25){\line(1,0){11}}\put(103,35){\line(1,0){11}}
\put(116,25){\line(1,0){11}}\put(116,35){\line(1,0){11}}
\put(82,23){\vector(0,1){6}}\put(84,23){\vector(0,1){6}} \put(86,23){\vector(0,1){6}}
\put(95,23){\vector(0,1){6}}\put(97,23){\vector(0,1){6}} \put(99,23){\vector(0,1){6}}
\put(108,23){\vector(0,1){6}}\put(110,23){\vector(0,1){6}} \put(112,23){\vector(0,1){6}}
\put(121,23){\vector(0,1){6}}\put(123,23){\vector(0,1){6}} \put(125,23){\vector(0,1){6}}
\put(82,33){\vector(0,1){6}}\put(85,33){\vector(0,1){6}}
\put(95,33){\vector(0,1){6}}\put(98,33){\vector(0,1){6}}
\put(108,33){\vector(0,1){6}}\put(111,33){\vector(0,1){6}}
\put(121,33){\vector(0,1){6}}\put(124,33){\vector(0,1){6}}
\put(0,16){$E_1$} \put(38,14){$E_2$}\put(51,14){$E_3$}\put(64,14){$E_4$}
\put(78,14){$E_5$}\put(91,14){$E_6$}\put(104,14){$E_7$}\put(117,14){$E_8$}
\put(75,36){$E_9$}\put(87,36){$E_{10}$} \put(100,36){$E_{11}$}\put(113,36){$E_{12}$}
\put(74,26){$E_{13}$}\put(87,26){$E_{14}$} \put(100,26){$E_{15}$}\put(113,26){$E_{16}$}
\end{picture}
\end{center}
We compute the numerical data.}
\\ ${}$ \begin{center} \begin{tabular}{|c|c|c|c|c|c|}
  \hline
  i=   & $1$ & $2,3,4$ & $5,6,7,8$ & $9,10,11,12$ & $13,14,15,16$  \\
\hline
  $N_i^P$ & 35 & 40 & 40 & 42 & 43 \\
  $N_i^Q$ & 35 & 35 & 35 & 35 & 35 \\
  $N_i$   & 0 & 5 & 5 & 7 & 8 \\
  $\nu_i$ & 2 & 3 & 3 & 4 & 4 \\
  \hline
\end{tabular}
\end{center}
\emph{We find that $Z_{top,f}(s)=-\frac{560s^3+1274s^2+767s+132}{2(7s+4)(s+1)(2s+1)}$ and so
the candidate pole $-3/5$ is not a pole.}\hfill $\square$
\end{example}
To prove the criterion for poles of plane curves, Veys points out in \cite{Veys} that there is an `ordered tree structure', i.e.\ the candidate poles are ordered increasingly along well-described paths in the tree. This structure fails in the meromorphic case. The other implication in Theorem \ref{thmVeys} still holds, i.e.\ if $s_0$ is a pole of $Z_{top,f}(s)$, then $s_0 = -\nu_i/N_i$ for
some exceptional curve $E_i$ intersecting at least three times other components, or
$s_0 = -1/N_j$ for some irreducible component $E_j$ of the strict transform of $f^{-1}\{0\}$. We prove this result now in the plane meromorphic case.

Let $\pi: X \rightarrow \mathcal{U}$ be an embedded resolution for $f=P/Q$ as in Section 1 and let $(E_i)_{i \in S}$ be the irreducible components of the total transform $\pi^{-1}(H)$.
Let $E(N,\nu)$ be an exceptional component of $\pi$ for which $N = N^P - N^Q \neq 0$ and suppose that $E$ intersects exactly $k$ other irreducible components of $\pi^{-1}(H)$, say $E_1(N_1,\nu_1),\ldots,$$E_k(N_k,\nu_k)$. For $i \in \{1,\cdots,k\}$, we set $\alpha_i=\nu_i-\frac{\nu}{N}N_i$.
\begin{lemma} \label{relations}
The following relations hold:
\begin{eqnarray}
\label{eqn1}
\sum_{i=1}^k N_i & =  & (-E \cdot E)N;  \\
\label{eqn2}
\sum_{i=1}^k \alpha_i & = & k-2.
\end{eqnarray}
\end{lemma}
\emph{Proof.}
We can follow the conceptual proof of Veys in \cite{Veys3} for holomorphic functions. In $Pic(X)$ we have
$\sum_{i \in S} N_i^P E_i=0$ and $\sum_{i \in S} N_i^Q E_i=0$. Hence $NE=-\sum{_l}(N_l^P-N_l^Q)E_l$, where $l$ runs over all components $E_l$, $l \in S$, except $E$. Intersecting with $E$ yields $N (E \cdot E) = -\sum_{i=1}^k N_i$.
This gives us Relation (1). \\
\noindent For Relation (2) we have that $K_X=\sum_{i \in S}(\nu_i-1)E_i$.
Hence in $Pic(X) \otimes \mathbb{Q}$ we have
\begin{eqnarray*}
K_X & = & \sum_{i \in S}(\nu_i-1)E_i - \frac{\nu}{N}(\sum_{i \in S} N_i E_i) \\
 & = & (\nu - 1) E + \sum_{l}(\nu_l-1)E_l - \frac{\nu}{N}(\sum_{l} N_l E_l) - \nu E \\
  & = & -E + \sum_l(\alpha_l-1)E_l
\end{eqnarray*}
where $l$ runs over all components $E_l$, $l \in S$, except $E$.

\noindent Intersecting with $E$ gives
$K_X \cdot E  =  - E \cdot E + \sum_{i=1}^k (\alpha_i-1)$
and by the adjunction formula we get
$-2 = \mbox{deg} K_E = \sum_{i=1}^k (\alpha_i-1)$.
\hfill $\blacksquare$

\

If $\nu/N \neq \nu_i/N_i$ for $i \in \{1,\cdots,k\}$, then the contribution $\mathcal{R}$ of the exceptional divisor $E$ to the residue of $Z_{top,f}(s)$  at $s=-\nu/N$ is given by
\[\mathcal{R}=\frac{1}{N}(2-k+\sum_{i=1}^k \frac{1}{\alpha_i}).\]

\noindent Completely analogous to the holomorphic case, one can now prove the following result.
\begin{proposition} \label{polethen}
Let $f=\frac{P}{Q}$ be a germ of a plane
meromorphic function. If $s_0$ is a pole of $Z_{top,f}(s)$, then
$s_0=-\nu_i/N_i$ for some exceptional component intersecting at
least three other components or $s_0=-1/N_i$ for some irreducible
component of the strict transform of $P$.
\end{proposition}
\emph{Proof.} If $s_0$ is not a candidate pole of order two, then
one can deduce from  Relation  (\ref{eqn2}) that the exceptional
components giving rise to the candidate pole $s_0$ and that have at most two intersections
do not contribute
to the residue. If $s_0$ is a candidate pole of order two, then
suppose $E_i$ intersects $E_j$ and $\nu_i/N_i=\nu_j/N_j$. If $E_i$
or $E_j$ is a component of the strict transform of $P$, then we are done.
Suppose now that $E_i$ and $E_j$ are exceptional components. If $E_i$ or $E_j$ has no other intersections,
then we get a contradiction with
Relation (\ref{eqn2}). If they intersect  exactly one other component, then
Relation (\ref{eqn2}) implies that this component also yields the
same candidate pole. By iterating the previous arguments, we can
conclude that there is a component of the strict transform or an
exceptional component intersecting at least three other components
that also gives rise to the candidate pole $-\nu_i/N_i$. \hfill
$\blacksquare$

\


\begin{center}
\textsc{3. Some facts about resolution for plane curves}
\end{center}
${}$\\
In this section we formulate some properties
of resolution of singularities of plane curves. These properties should be well known but we recall them here because they are important to understand well the resolution of plane meromorphic germs.
This section will serve as a preparation for Section 4 where we will treat the monodromy conjecture for plane meromorphic germs.
\\ \\ \noindent
\textbf{(1)} \textbf{Notation.} Let us fix some terminology about the dual resolution graph $\mathcal{G}$ (of plane curves or of plane meromorphic germs). For a vertex of $\mathcal{G}$, let its \emph{valence} be the number of irreducible components that meet the irreducible component corresponding to that vertex. A \emph{bamboo} leaving out of a vertex of $\mathcal{G}$ is a connected path in $\mathcal{G}$ starting in that vertex, having an end point with valence 1 and in which the vertices in between have all valence 2.
A bamboo is said to be a \emph{primitive bamboo} if its valence one vertex does not correspond to an irreducible component
of the strict transform (i.e.\ it is not an arrowhead). The \emph{branches} leaving out of a vertex of $\mathcal{G}$
are the connected components of $\mathcal{G}$ minus that vertex. Such a branch is called
primitive if none of its valence one vertices corresponds to an irreducible component of the strict transform
(i.e.\ the branch does not contain arrowheads).
\\ \\ \noindent
\textbf{(2)} \textbf{Dual resolution graphs.} The following properties of dual resolution graphs of plane curve singularities will be useful. See \cite[III. Section 8.4, Prop. 16]{brkn} or \cite[Section 5.4, Theorem 5.4.5]{jongpfister} for a proof.
\begin{proposition} \label{remminimal1}
For the dual resolution graph $\mathcal{G}$ of the minimal embedded resolution of a plane curve the following statements hold:
\begin{enumerate}
\item The graph $\mathcal{G}$ is finite and does not contain cycles;
\item For every vertex of $\mathcal{G}$ there is at most one primitive bamboo leaving out of that vertex, except for possibly one vertex from which at most two primitive bamboos leave;
\item For every vertex of $\mathcal{G}$ there is at most one primitive branch leaving out of that vertex which is not a bamboo.
\end{enumerate}
\end{proposition}

\begin{remark} \label{remminimal2}
\emph{To create the minimal embedded resolution of a meromorphic germ $P/Q$, notice that one can take the minimal embedded resolution of $PQ$ to which one adds the necessary blow ups to create the dicritical components. Notice that the blow ups (the intersections of irreducible components of $\pi^{-1}\{PQ=0\}$) necessary to create the dicritical components only create vertices with valence 2. This implies that the statements in Proposition \ref{remminimal1} still hold for the dual minimal embedded resolution graph of a plane meromorphic germ.}
\end{remark}

\noindent \textbf{(3) } \textbf{Properties of the resolution data.} Let $f:(\C^2,0) \rightarrow (\C,0)$ be a holomorphic germ. Let $E(N,\nu)$ be an exceptional component with $N \neq 0$ and suppose that $E$ intersects exactly $k$ other irreducible components of $(f \circ \pi)^{-1}\{0\}$, say $E_1(N_1,\nu_1),\ldots,$$E_k(N_k,\nu_k)$. For $i \in \{1,\cdots,k\}$, we set $\alpha_i=\nu_i-\frac{\nu}{N}N_i$.
The bounds on the numbers $\alpha_i$ play an essential role in the proof by Rodrigues of the monodromy conjecture for plane curves (see \cite{bart}).
In the holomorphic case Loeser stated the following bounds:
\begin{proposition} {\rm \cite[Proposition II.3.1]{loeser1}} \label{proploeser}
Suppose that $E$ is an exceptional component in the minimal embedded resolution of a plane curve which intersects at least two irreducible components in the total transform, one of them being $E_j$. Then $-1 < \alpha_{j}< 1$. If the embedded resolution of the plane curve is not minimal, then $-1 \leq \alpha_{j} \leq 1$.
\end{proposition}
\begin{remark}
\emph{If the exceptional component $E$ intersects exactly one component, say $E_j$, then it follows from Relation (2) for holomorphic germs (see \cite[Lemma 2.3]{Veys2}) that $\alpha_j=-1$.}
\end{remark}
In our proof of the monodromy conjecture for plane meromorphic germs we will follow the same strategy as Rodrigues but we will have to study carefully the numbers $\alpha_i$. Indeed, in the case of meromorphic germs in general the numbers $\alpha_i$ are not bounded even in the minimal embedded resolution (see Example 1 with $E=E_8$ and $E_j=E_1$, then $\alpha_1=-5/2$). In what follows we will look at subgraphs of the dual resolution graph of a plane meromorphic function in which the numbers $\alpha_i$ are contained in the interval $(-1,1)$.

Let $\pi : X \rightarrow \mathcal{U}$ be an embedded resolution of the germ of the plane meromorphic function $f = \frac{P}{Q}$ at 0 and let $(E_i)_{i \in S}$ be the irreducible components of the total transform of $f$. Denote by $\mathcal{E}$ the  exceptional divisor of $\pi$; i.e. $\mathcal{E}= \cup_{i\in S}E_i \cap \pi^{-1}(0)$, and by $\mathcal{D}$ the union of the dicritical components of $\pi$.

In order to prove the monodromy conjecture for plane meromorphic germs we need to take into account some features of  the connected components of the closure  $(\mathcal{E} - \mathcal{D})^c$ of $\mathcal{E} - \mathcal{D}$. These components where studied by Delgado and Maugendre  in relation to the geometry of pencils of plane curve singularities (see \cite{delgadom}). Here we recall some properties needed in our proof.

Let $E(N,\nu)$ be an exceptional component of $\pi$ with $N = N^P - N^Q \neq 0$ and suppose that $E$ intersects exactly $k$ other irreducible components of $(f \circ \pi)^{-1}\{0\}$, say $E_1(N_1,\nu_1),\ldots,$$E_k(N_k,\nu_k)$.
Recall that $\pi$ is also a resolution of the curves $P$ and $Q$ and that the corresponding resolution data are related by $N_i = N_i^P-N_i^Q$.  From the identity  \cite[2.2 (ii)]{Veys3}, which is the version of  (\ref{eqn2}) for resolutions of plane curve singularities, we have that
$\sum_{i=1}^k \alpha^P_i = k-2 = \sum_{i=1}^k \alpha^Q_k.$ Equivalently, we have the identity
\begin{equation}\label{eqn3}\sum_{i=1}^kN_i^P/N^P = \sum_{i=1}^k N_i^Q/N^Q,\end{equation}
which implies the following result.

\begin{proposition}\label{del1} {\rm \cite[Proposition 1]{delgadom}} Assume that none of the components $E_1,\ldots, E_k$ belongs to the strict transform of $ \{P_{|_\mathcal{U}}\} \cup \{Q_{|_\mathcal{U}}\}$. Then,  there exists an index $r \in \{1, \dots, k\}$ such that $N_r^P N^Q > N_r^Q N^P$ if and only if there exists an index $s \in \{1, \dots, k\}$ such that $N_s^P N^Q < N_s^Q N^P$.\end{proposition}

\begin{corollary}\label{remminstrict}{\rm \cite[Corollary 1]{delgadom}}\emph{Proposition \ref{del1} implies that the quotient $N^P/N^Q$ is constant on a primitive bamboo.}
\end{corollary}
\begin{example}\emph{
Let us consider two irreducible plane curves $P(x,y)=x^{103}-y^{24}=0$ and $Q(x,y)=x^{30}-y^7=0$.
The dual resolution graph of $P/Q$ is as follows
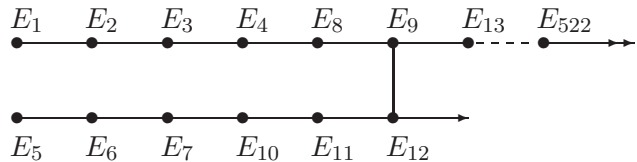
\begin{center}
\begin{picture}(90,30)(0,0)
\put(9,19){$\bullet$} \put(19,19){$\bullet$} \put(29,19){$\bullet$} \put(39,19){$\bullet$} \put(49,19){$\bullet$} \put(59,19){$\bullet$}
\put(60,10){\vector(1,0){10}}\put(10,20){\line(1,0){50}}
\drawline(60,20)(60,10)
\put(9,22){$E_{1}$}\put(19,22){$E_{2}$}\put(29,22){$E_{3}$}\put(39,22){$E_{4}$}\put(49,22){$E_{8}$}
\put(59,22){$E_{9}$}
\put(80,20){\vector(1,0){10}}\put(80,20){\vector(1,0){12}}
\put(9,9){$\bullet$} \put(19,9){$\bullet$} \put(29,9){$\bullet$} \put(39,9){$\bullet$}
\put(49,9){$\bullet$} \put(59,9){$\bullet$} \put(69,19){$\bullet$} \drawline(59,20)(69,20)
\drawline(69,20)(70,20)\drawline(71,20)(72,20)\drawline(73,20)(74,20)\drawline(75,20)(76,20)\drawline(77,20)(78,20)
\put(79,19){$\bullet$}
\put(10,10){\line(1,0){50}}
\put(9,5){$E_{5}$}\put(19,5){$E_{6}$}\put(29,5){$E_{7}$}\put(39,5){$E_{10}$}\put(49,5){$E_{11}$}
\put(59,5){$E_{12}$}\put(69,22){$E_{13}$}\put(79,22){$E_{522}$}
\end{picture}
\end{center}
where $E_{13},\ldots,E_{522}$ are exceptional components created to get a dicritical component.
The numerical data for the components $E_1,\ldots,E_{12}$ are given in the following table
$$\begin{array}{|c|c|c|c|c|c|c|c|c|c|c|c|c|}
\hline
&E_1&E_2&E_3&E_4&E_5&E_6&E_7&E_8&E_9&E_{10}&E_{11}&E_{12}\\
\hline
N^P_i&24&48&72&96&103&206&309&408&720&1030&1751&2472\\
\hline
N^Q_i&7&14&21&28&30&60&90&119&210&300&510&720\\
\hline
\nu_i&2&3&4&5&6&11&16&21&37&53&90&127\\
\hline
\end{array}$$
Corollary \ref{remminstrict} tells us that
\[\frac{24}{7}= \frac{N_1^P}{N_1^Q}= \frac{N_2^P}{N_2^Q}=\frac{N_3^P}{N_3^Q}=\frac{N_4^P}{N_4^Q}=\frac{N_8^P}{N_8^Q}=\frac{N_9^P}{N_9^Q} \qquad \mbox{ and} \]
\[\frac{103}{30}= \frac{N_5^P}{N_5^Q}= \frac{N_6^P}{N_6^Q}=\frac{N_7^P}{N_7^Q}=\frac{N_{10}^P}{N_{10}^Q}=\frac{N_{11}^P}{N_{11}^Q}=\frac{N_{12}^P}{N_{12}^Q}.\]}
\hfill $\square$
\end{example}

\

\noindent \textbf{(4)} \textbf{Some special subgraphs.}  If $\mathcal{K}$ is a connected component of $(\mathcal{E} - \mathcal{D})^c$, then $f \circ \pi$ is constant along $\mathcal{K}$. Let us assume that $(f \circ \pi)_{|_\mathcal{K}}$ is constantly $0$.
\begin{proposition}\label{del2} {\rm \cite[Theorem 1]{delgadom}} $\mathcal{K}$ intersects the strict transform of the hypersurface  $\{P_{|_\mathcal{U}}=0\}$.
\end{proposition}
\emph{Proof.} If there are no dicritical components, we are in the case of plane curves and obviously some exceptional component intersects a component of the strict transform of $\{P_{|_\mathcal{U}}=0\}$. So suppose that the connected component $\mathcal{K}$ meets a dicritical component.
Since $\mathcal{K}$ is finite we can take an element $E \in \mathcal{K}$ such that $\frac{N^P}{N^Q}$ is maximal in $\mathcal{K}$.  If $E$ does not intersect the strict transform of the hypersurface  $\{P_{|_\mathcal{U}}=0\}$ and if $E_1,\ldots, E_k$ are the components that intersect $E$, then Proposition \ref{del1} implies that $\frac{N^P}{N^Q}=\frac{N^P_1}{N^Q_1}= \cdots=\frac{N^P_r}{N^Q_r}$. In this way, one gets that the quotient
$\frac{N^P}{N^Q}$ is constant in all $\mathcal{K}$, provided that $\mathcal{K}$ does not intersect the strict transform of the hypersurface  $\{P_{|_\mathcal{U}}=0\}$. However, again Proposition \ref{del1} prevents that any of the components in $\mathcal{K}$ intersects a dicritical component because $\frac{N^P}{N^Q} > 1$ for all components in $\mathcal{K}$ as $(f \circ \pi)_{|_\mathcal{K}}$ is constantly $0$ and $\frac{N^P}{N^Q} = 1$ for all dicritical components. Thus we reach a contradiction.
\hfill
$\blacksquare$ ${}$ \\ ${}$

\noindent From now on we take $\pi : X \rightarrow \mathcal{U}$ the minimal embedded resolution of the plane meromorphic function $f = \frac{P}{Q}$.
\begin{proposition} \label{prophelp} Let $P,Q: (\mathbb{C}^2,0) \rightarrow (\mathbb{C},0)$ and let $\mathcal{G}$ be the dual resolution graph corresponding to the minimal
embedded resolution $\pi$ of $P/Q$. Let $\mathcal{P}$ be a connected subgraph of $\mathcal{G}$ that contains exactly one vertex which is connected to vertices
not belonging to $\mathcal{P}$, 
that does not contain vertices corresponding to components of the strict transform and on which $(f \circ \pi)_{|_\mathcal{P}}$ is constantly $0$.
Let $E_i$ and $E_j$ be two exceptional components which are represented by adjacent vertices in $\mathcal{P}$. If $E_i$ is intersected by at least one more component, then $-1 < \alpha_{i,j} := \nu_j - \frac{\nu_i}{N_i} N_j < 1$.
\end{proposition}
\emph{Proof.} Let us first argue on how the subgraph $\mathcal{P}$ then looks like.
Proposition \ref{remminimal1} and Remark \ref{remminimal2} imply that $\mathcal{P}$ has the following shape, where we only draw the vertices with valence different from two.

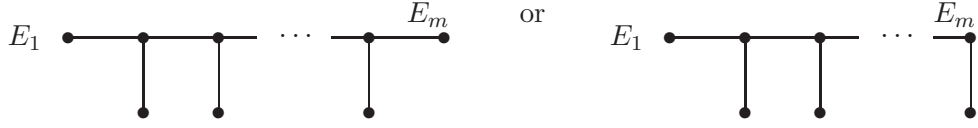
\begin{center}
\begin{picture}(120,25)(0,0)
\put(9,19){$\bullet$} \put(19,19){$\bullet$}  \put(39,19){$\bullet$} \put(49,19){$\bullet$}
\put(19,9){$\bullet$}
\put(-1,19){$\bullet$} \put(9,9){$\bullet$}
\put(39,9){$\bullet$}

\put(0,20){\line(1,0){25}}
\put(28,20){\dots}
\put(35,20){\line(1,0){15}}
\put(10,10){\line(0,1){10}}
\put(20,10){\line(0,1){10}}
\put(40,10){\line(0,1){10}}


\put(-8,19){$E_{1}$}
\put(45,22){$E_{m}$}

\put(60,22){or}

\put(89,19){$\bullet$} \put(99,19){$\bullet$}  \put(119,19){$\bullet$}
\put(99,9){$\bullet$}
\put(79,19){$\bullet$} \put(89,9){$\bullet$}
\put(119,9){$\bullet$}

\put(80,20){\line(1,0){25}}
\put(108,20){\dots}
\put(115,20){\line(1,0){5}}
\put(90,10){\line(0,1){10}}
\put(100,10){\line(0,1){10}}
\put(120,10){\line(0,1){10}}

\put(72,19){$E_{1}$}
\put(115,22){$E_{m}$}

\end{picture}
\end{center}
\vspace*{-0.8cm} with particular cases
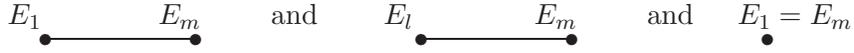
\begin{center}
\begin{picture}(90,20)(15,10)
\put(9,19){$\bullet$}  \put(29,19){$\bullet$}
\put(10,20){\line(1,0){20}}
\put(5,22){$E_{1}$}
\put(25,22){$E_{m}$}

\put (40,22){and}

\put(59,19){$\bullet$}  \put(79,19){$\bullet$}
\put(60,20){\line(1,0){20}}
\put(55,22){$E_{l}$}
\put(75,22){$E_{m}$}

\put (90,22){and}
\put(105,19){$\bullet$} \put(102,22){$E_1=E_{m}$}

\end{picture}
\end{center}
\vspace*{-0.8cm} Here $E_1$ denotes the exceptional component corresponding to the first blow up of $\pi$, the component $E_m$ is the component that intersects some components outside of $\mathcal{P}$ and $E_l \neq E_1$. In the last case where $\mathcal{P}$ only consists of one exceptional divisor, $E_m=E_1$.

 Next notice that $\mathcal{P}$ was already a subgraph of the minimal resolution of $PQ$. Indeed, suppose that any vertex $E$ in $\mathcal{P}$ is created to pass from the minimal embedded resolution of $PQ$ to the embedded resolution of $f=\frac{P}{Q}$. According to Remark \ref{remminimal2}, $E$ has valence 2 and erasing $E$ divides the resolution graph into two disjoint parts, one of then containing at least irreducible components of the strict transform of  $\pi^{-1}\{P=0\}$ and the other containing at least irreducible components of the strict transform $\pi^{-1}\{Q=0\}$. This is a contradiction as we supposed that $\mathcal{P}$  contains no irreducible components of the strict transform of  $\pi^{-1}\{P=0\}$ neither of $\pi^{-1}\{Q=0\}$.

 Therefore we can assume that $E_i$ and $E_j$ appeared in the minimal resolution graph of $PQ$ and Proposition \ref{proploeser} tells us then that
\begin{eqnarray}\label{ineq}
-1 < \nu_j - \frac{\nu_i}{N_i^P+N_i^Q} (N_j^P+N_j^Q )< 1.
\end{eqnarray}

Taking into account the form of the subgraph $\mathcal{P}$, we can now prove that $N^P/N^Q$ is constant on the whole graph $\mathcal{P}$.
Corollary \ref{remminstrict} deals with the case that $\mathcal{P}$ is a primitive bamboo. Otherwise let  $E$ be an exceptional component in $\mathcal{P}$ corresponding to a vertex of valence equal to $3$ from which two bamboos are leaving. Then it follows from Corollary \ref{remminstrict} that the quotient $N_i^P/N_i^Q$ is constant on these bamboos and equal to $N^P/N^Q$.  Denote by  $E_1,E_2$ and $E_3$ the exceptional components intersecting $E$.  Since $N^P/N^Q=N_1^P/N_1^Q=N_2^P/N_2^Q$, we can deduce that
$N^P/N^Q=N_3^P/N_3^Q$. Indeed, Relation (\ref{eqn2}) in the holomorphic case gives
\[ \nu_1 -\frac{\nu}{N^P}N_1^P + \nu_2 -\frac{\nu}{N^P}N_2^P + \nu_3 -\frac{\nu}{N^P}N_3^P= 1 = \nu_1 -\frac{\nu}{N^Q}N_1^Q + \nu_2 -\frac{\nu}{N^Q}N_2^Q + \nu_3 -\frac{\nu}{N^Q}N_3^Q,\] from which we can deduce $N^P/N^Q=N_3^P/N_3^Q$. Next we can propagate this equality along a chain of vertices of valence 2.  Let $E_j$ be a vertex of has valence 2 intersecting
$E_i$ and $E_k=E$. Then we have that
\[ \nu_i -\frac{\nu_j}{N_j^P}N_i^P + \nu_k -\frac{\nu_j}{N_j^P}N_k^P = 0 = \nu_i -\frac{\nu_j}{N_j^Q}N_i^Q + \nu_k -\frac{\nu_j}{N_j^Q}N_k^Q.\] Now we can deduce that $N_i^P/N_i^Q = N_j^P/N_j^Q$.
We can iterate this until we arrive at another vertex of valence equal to 3, from which a primitive bamboo is leaving. Using again the previous argument and taking into account the shape of $\mathcal{P}$, we deduce that the quotients $N_i^P/N_i^Q$ are equal for all $E_i$ in $\mathcal{P}$.

Finally, since $\frac{N_i^{P}}{N_j^P}=\frac{N_i^Q}{N_j^Q}$ and hence
$\frac{N_i^P+N_i^Q}{N_j^P+N_j^Q}=\frac{N_i^{P}}{N_j^P}=\frac{N_i}{N_j}$. We substitute in (\ref{ineq}) and get
$-1 < \nu_j - \frac{\nu_i}{N_i}N_j < 1$.
\hfill
$\blacksquare$ ${}$ \\ ${}$

\noindent \textbf{(5)} \textbf{Topology of trees of nonsingular rational curves}  The following lemma was shown by Rodrigues to give his own proof of the monodromy conjecture for plane curves.
\begin{lemma} {\rm \cite[Lemma(2.2)]{bart}}\label{bart}
Let $ \cup_{i=1}^r E_i$ be a tree of nonsingular rational curves on a nonsingular surface, then
\[\sum_{i=1}^r \chi(E_i^{\circ})=2.\]
\end{lemma}

\begin{center}
\textsc{4. The monodromy conjecture for plane meromorphic germs}
\end{center}
${}$\\
By analogy with the holomorphic case, we propose the monodromy conjecture for meromorphic germs as follows.

\begin{conjecture} (Monodromy Conjecture for meromorphic germs)
Let $f=\frac{P}{Q}$ be a germ of a meromorphic function on $(\C^{n+1},0)$ and let $s_0$ be a pole of $Z_{top,f}$, then $e^{2\pi i s_0}$ is an eigenvalue of the monodromy transformation $h_{f}$ at some point of $P^{-1}\{0\}$ in a neighbourhood of $0$.
\end{conjecture}
\begin{theorem}
Let $f=\frac{P}{Q}$ be a germ of a meromorphic function on $(\C^{2},0)$ and let $s_0$ be a pole of $Z_{top,f}(s)$, then $e^{2 \pi i s_0}$ is an eigenvalue of the monodromy transformation $h_{f}$ at some point of $P^{-1}\{0\}$ in a neighbourhood of $0$.
\end{theorem}
\emph{Proof.} Let $\pi: X \rightarrow \mathcal{U}$ be the minimal embedded resolution of $f$. By
Proposition \ref{polethen} it follows that either $s_0=-1/N$ for
some irreducible component $E(N,1)$ of the strict transform of $P$ or $s_0=-\nu_i/N_i$ for some
exceptional divisor $E_i$ of $\pi$ with $\chi(E_i^{\circ}) < 0$.
In the first case we find the eigenvalue $e^{2\pi i s_0}$ at any point of $\{P=0\}$ in a
punctured neighbourhood of $0$. So suppose now that $s_0=-\nu_i/N_i$ for some
exceptional divisor $E_i$ of $\pi$ with $\chi(E_i^{\circ}) < 0$. We write $\nu_i/N_i=a/d$ with gcd$(a,d)=1$. We can suppose that there is no component $E_S(N_S,\nu_S)$ of the strict transform of $P$ for which $d \mid N_S$. Otherwise  $e^{2 \pi i s_0}$ is an eigenvalue of monodromy at any point of $\{P=0\}$ in a punctured neighbourhood of $0$.

Let us define
\[C_d:=\{E_k \textnormal{ exceptional divisor of } \pi  \textnormal{ such that } N_k > 0 \textnormal{  and } d \mid N_k\}.\]
 Now we study $\sum_{E_k \in \mathcal{K}} \chi(E_k^\circ)$ for each  connected component $\mathcal{K}$ of $C_d$. Firstly,  Proposition \ref{del2} shows that there is no such component $\mathcal{K}$  intersecting only dicritical divisors.
Notice that an exceptional component $E$ with $N=0$ intersecting such a connected component $\mathcal{K}$ is also dicritical because $(f \circ \pi)_{|_{E}}$ then reaches at least two different values.
Hence, each connected component $\mathcal{K}$ in $C_d$ does contain a component $E_l$ that intersects some  irreducible component $E$ of the total transform $\pi^{-1}(H)$ outside $\mathcal{K}$ with $N \not =0$.
Relation (1) in Lemma \ref{relations} implies that there is at least another irreducible component of $\pi^{-1}(H)$ outside $\mathcal{K}$ satisfying these properties. Lemma \ref{bart} then implies that $\sum_{E_k \in
\mathcal{K}}\chi(E_k^{\circ}) \leq 0$ for each connected component $\mathcal{K}$ of $C_d$.
Now the result follows if for at least one such component  $\sum_{E_k \in
\mathcal{K}}\chi(E_k^{\circ}) < 0$ holds. We show that this is the case for the component $\mathcal{K}$ of $C_d$ that contains the exceptional component $E_i$. We denote this component by $\tilde{\mathcal{N}}$. Let us argue by contradiction.  We suppose that $\sum_{E_k \in
\tilde{\mathcal{N}}}\chi(E_k^{\circ})$ is equal to $0$. This assumption implies that $\tilde{\mathcal{N}}$ contains  only one exceptional component intersecting irreducible components of $\pi^{-1}(H)$ not belonging to $\tilde{\mathcal{N}}$ and with $N \not =0$. Furthermore the number of the intersections must be exactly two. Thus  $\tilde{\mathcal{N}}$ satisfies the conditions of the component $\mathcal{P}$ in Proposition \ref{prophelp}. In particular, we can disregard the case of $\tilde{\mathcal{N}}$ being a singleton (i.e. $\tilde{\mathcal{N}}$ only consists of $E_i$) because $\chi(E_i^{\circ}) < 0$. Remark that in this case $e^{-2 \pi ia/d }$ is an eigenvalue of monodromy at the origin.
If $\tilde{\mathcal{N}}$ is not a singleton, then let $E_j$ be an exceptional component in $\tilde{\mathcal{N}}$ intersecting $E_i$. As $\alpha_{i,j}=\nu_j-\frac{\nu_i}{N_i}N_j= \nu_j-\frac{a}{d}N_j \in \mathbb{Z}$, it follows by Proposition \ref{prophelp} that $\alpha_{i,j}=0$. Hence $\nu_j/N_j=a/d$. On the other hand, as $\alpha_{i,j}=0$, it follows from Relation (2) that $E_j$ intersects exactly one other component $E_k$. If $E_k \in \tilde{\mathcal{N}}$, then analogously we find $\nu_k/N_k=a/d$ and using Relation (2) we see that $E_k$ intersects exactly one other component. As the graph $\mathcal{G}$ is finite, we must end with an exceptional component not in $\tilde{\mathcal{N}}$.
However, as $\tilde{\mathcal{N}}$ only contains one exceptional component which intersects exceptional components not in $\tilde{\mathcal{N}}$ and as $E_i$ is intersecting at least three components,
we can do the same reasoning for the other exceptional components intersecting $E_i$ and hence get a contradiction.
\hfill $\blacksquare$

${}$\\ \\
\noindent Manuel Gonz\'alez Villa; Universit\"{a}t Heidelberg, MATCH, Im Neuenheimerfeld 288,
69120 Heidelberg, email: \texttt{villa@mathi.uni-heidelberg.de}
\\ \\
Ann Lemahieu; Universit\'e Lille 1, U.F.R. de Math\'ematiques, 59655 Villeneuve d'Ascq C\'edex, email: \texttt{ann.lemahieu@math.univ-lille1.fr}
\\ \\
Acknowledgement:  The authors thank Antonio Campillo, Alejandro Melle and Wim Veys for their helpful comments.
The research was partially supported by the Fund of
Scientific Research - Flanders, MEC PN I+D+I MTM2007-64704, MCI-Spain grant MTM2010-21740-C02 and the ANR `SUSI' project (ANR-12-JS01-0002-01). The second author is very grateful to the Universidad Complutense de Madrid (Spain) for its hospitality.


\begin{thebibliography}{00}

\bibitem[A]{arnold} V.I. Arnol'd, \emph{Singularities of fractions and the behavior of polynomials at infinity},
(Russian) Tr. Mat. Inst. Steklova \textbf{221} (1998), 48--68; translation in
Proc. Steklov Inst. Math. 1998, no. \textbf{2} (221), 40--59.

\bibitem[ACLM1]{luengo} E. Artal Bartolo, P. Cassou-Nogu\`{e}s, I.
Luengo, A. Melle Hern\'{a}ndez, \emph{Monodromy conjecture for some
surface singularities}, Ann. Scient. Ec. Norm. Sup. \textbf{35}
(2002), 605--640.

\bibitem[ACLM2]{luengo2} E. Artal Bartolo, P. Cassou-Nogu\`{e}s, I.
Luengo, A. Melle Hern\'{a}ndez, \emph{Quasi-ordinary power series
and their zeta functions}, Memoirs of the A.M.S. \textbf{178}, no.
841 (2005).

\bibitem[BK]{brkn} E. Brieskorn and H. Kn\"{o}rrer, \emph{Plane algebraic curves}, Birkh\"{a}user Verlag Basel, 1986.

\bibitem[DM]{delgadom} F. Delgado and H. Maugendre, \emph{Special Fibres and Critical Locus for a Pencil of Plane Curve Singularities}, Comp. Math. \textbf{136} (2003) 69--87.


\bibitem[dJP]{jongpfister} T. de Jong and G. Pfister, \emph{Local analytic geometry}, Advanced Lectures in Mathematics, Vieweg, 2000.

\bibitem[GZLM1]{meromorphic} S.M. Gusein-Zade, I. Luengo and A. Melle-Hern\'andez, \emph{Zeta functions for germs of meromorphic functions, and Newton diagrams}, Functional Analysis and its applications, Vol. \textbf{32}, No. 2, 1998.

\bibitem[GZLM2]{meromorphic2} S.M. Gusein-Zade, I. Luengo and A. Melle-Hern\'andez, \emph{On the topology of germs of meromorphic functions and its applications}, (Russian) Algebra i Analiz \textbf{11} (1999), no. 5, 92--99; translation in St. Petersburg Math. J. \textbf{11} (2000), no. 5, 775--780.

\bibitem[LVa]{liseann} A. Lemahieu and L. Van Proeyen, \emph{Monodromy conjecture for nondegenerate surface singularities}, Trans. of the A.M.S. \textbf{363} (2011), 4801--4829.

\bibitem[LV]{monodconjtoric} A. Lemahieu and W. Veys, \emph{Zeta functions and monodromy for surfaces that are general for a toric idealistic
cluster}, Int. Math. Res. Notices, ID rnn 122 (2009), 52 pages.

\bibitem[L1]{loeser1} F. Loeser, \emph{Fonctions d'Igusa p-adiques
et polyn\^{o}mes de Bernstein.} Amer. J. Math. \textbf{110} (1988),
1--22.

\bibitem[L2]{loeser2} F. Loeser, \emph{Fonctions d'Igusa p-adiques, polyn\^{o}mes de Bernstein,
 et poly\`{e}dres de Newton}, J. reine angew. Math. \textbf{412} (1990),
 75--96.

\bibitem[NV1]{wimnemethi1} A. N\'emethi and W. Veys, \emph{Monodromy eigenvalues are induced by
poles of zeta functions: the irreducible curve case}, Bull. Lond. Math.
Soc. \textbf{42} (2010), no. 2, 312--322.

\bibitem[NV2]{wimnemethi2} A. N\'emethi and W. Veys, \emph{Generalized monodromy conjecture in dimension two}, Geometry and Topology \textbf{16} (2012), 155--217.


\bibitem[R1]{raibaut} M. Raibaut, \emph{Motivic Milnor fibers of a rational function},  C. R. Math. Acad. Sci. Paris {\bf 350} (2012), no. 9-10, 519--524.

\bibitem[R2]{raibautb} M. Raibaut, \emph{Motivic Milnor fibers of a rational function},  To appear at Revista Matem\'atica Complutense.

\bibitem[Ro]{bart} B. Rodrigues, \emph{On the monodromy conjecture for curves on normal surfaces}, Math. Proc. of the Cambridge Philosophical Society \textbf{136} (2004) 313--324.


\bibitem[Ve1]{Veys} W. Veys, \emph{Determination of the poles of the topological zeta function for curves},
Manuscripta Math. \textbf{87} (1995), 435--448.

\bibitem[Ve2]{Veys2} W. Veys, \emph{The topological zeta function associated to a function on a normal surface germ}, Topology \textbf{38} (1999), 439--456.

\bibitem[Ve3]{Veys3} W. Veys, \emph{Embedded resolution of singularities and Igusa's local zeta function}, Academiae Analecta: Mededelingen van de Koninklijke Academie voor Wetenschappen, Letteren, (2001), 1--56.

\bibitem[Ve4]{Veys4} W. Veys, \emph{Monodromy eigenvalues and zeta functions with
differential forms}, Adv. Math. \textbf{213} (2007), no. 1, 341--357.


\end{thebibliography}
\end{document}